\documentclass[leqno,12pt]{amsart}
\usepackage{amssymb} 
\theoremstyle{plain} 
\newtheorem{theorem}{\indent\sc Theorem}[section] 

\newtheorem{proposition}[theorem]{\indent\sc Proposition}

\theoremstyle{definition} 

\begin{document}

\title[Coupled Hamiltonian systems]{Coupled Hamiltonian systems with extended affine Weyl group symmetry of type $D_3^{(2)}$ \\}
\author{Yusuke Sasano }

\renewcommand{\thefootnote}{\fnsymbol{footnote}}
\footnote[0]{2000\textit{ Mathematics Subjet Classification}.
34M55; 34M45; 58F05; 32S65.}

\maketitle

\begin{abstract}
We find a two-parameter family of ordinary differential systems in dimension five with the affine Weyl group symmetry of type $D_3^{(2)}$. We show its symmetry and holomorphy conditions. This is the second example which gave higher order Painlev\'e type systems of type $D_{3}^{(2)}$. By obtaining its first integrals of polynomial type, we can obtain a two-parameter family of coupled Hamiltonian systems in dimension four with the polynomial Hamiltonian.
\end{abstract}

\section{Introduction}

In this paper, we find a 2-parameter family of ordinary differential systems in dimension five with the affine Weyl group symmetry of type $D_3^{(2)}$ explicitly given by
\begin{equation}\label{eq:11}
  \left\{
  \begin{aligned}
\frac{dx}{dt}=&-(xw-\alpha_2)x+\frac{1}{2},\\
\frac{dy}{dt}=&(xw+zq-1)y+\alpha_1 wq,\\
\frac{dz}{dt}=&-(zq-\alpha_0)z-\frac{\eta}{2},\\
\frac{dw}{dt}=&(xw-zq-\alpha_2)w+yz,\\
\frac{dq}{dt}=&(zq-xw-\alpha_0)q+xy.
   \end{aligned}
  \right. 
\end{equation}
Here $x,y,z,w$ and $q$ denote unknown complex variables, and $\alpha_0,\alpha_1,\alpha_2$ are complex parameters satisfying the relation:
\begin{equation}\label{eq:12}
\alpha_0+\alpha_1+\alpha_2=1.
\end{equation}

This is the second example which gave higher order Painlev\'e type systems of type $D_{3}^{(2)}$.

We also remark that 2-coupled Painlev\'e II system in dimension four given in the paper \cite{Sasano3} admits the affine Weyl group symmetry of type $D_3^{(2)}$ as the group of its B{\"a}cklund transformations, whose generators $s_0,s_1,s_2$ are determined by the invariant divisors.

On the other hand, the B{\"a}cklund transformations $s_0,s_2$ of this system do not satisfy so except for the transformation $s_1$ (see Theorem \ref{th:1}).

We show its symmetry and holomorphy conditions. By obtaining its first integrals of polynomial type, we can obtain a two-parameter family of coupled Hamiltonian systems in dimension four with the polynomial Hamiltonian. Each principal part of the Hamiltonian $H$ has its first integral, respectively. Nevertheless, the Hamiltonian $H$ itself is not its first integral.

\section{Symmetry and holomorphy conditions}

\begin{theorem}
Let us consider the following ordinary differential system in the polynomial class\rm{:\rm}
\begin{align*}
&\frac{dx}{dt}=f_1(x,y,z,w,q), \cdots, \frac{dq}{dt}=f_5(x,y,z,w,q) \quad (f_i \in {\Bbb C}(t)[x,y,z,w,q]).
\end{align*}
We assume that

$(A1)$ $deg(f_i)=3$ with respect to $x,y,z,w,q$.

$(A2)$ The right-hand side of this system becomes again a polynomial in each coordinate system $(x_i,y_i,z_i,w_i,q_i) \ (i=0,1,2):$

\begin{align}
\begin{split}
0) \ &x_0=x, \quad y_0=y-\frac{2\alpha_0 w}{z}+\frac{\eta w}{z^2} \quad z_0=z, \quad w_0=w, \quad q_0=q-\frac{2\alpha_0}{z}+\frac{\eta}{z^2},\\
1) \ &x_1=x+\frac{\alpha_1 q}{y}, \quad y_1=y, \quad z_1=z+\frac{\alpha_1 w}{y}, \ w_1=w, \quad q_1=q,\\
2) \ &x_2=x, \quad y_2=y-\frac{2\alpha_2 q}{x}-\frac{q}{x^2} \quad z_2=z, \quad w_2=w-\frac{2\alpha_0}{x}-\frac{1}{x^2}, \quad q_2=q.
\end{split}
\end{align}
Then such a system coincides with the system \eqref{eq:11}.
\end{theorem}
We note that these transition functions satisfy the condition{\rm:\rm}
\begin{align*}
&dx_i \wedge dy_i \wedge dz_i \wedge dw_i \wedge dq_i=dx \wedge dy \wedge dz \wedge dw \wedge dq \quad (i=0,1,2).
\end{align*}

\begin{theorem}\label{th:1}
The system \eqref{eq:11} admits the affine Weyl group symmetry of type $D_3^{(2)}$ as the group of its B{\"a}cklund transformations, whose generators $s_0,s_1,s_2$ defined as follows$:$ with {\it the notation} $(*):=(x,y,z,w,q,\eta;\alpha_0,\alpha_1,\alpha_2)$,
\begin{align*}
s_0:(*) \rightarrow &\left(x,y-\frac{2\alpha_0 w}{z}+\frac{\eta w}{z^2},z,w,q-\frac{2\alpha_0}{z}+\frac{\eta}{z^2},-\eta;-\alpha_0,\alpha_1+2\alpha_0,\alpha_2 \right),\\
s_1:(*) \rightarrow &\left(x+\frac{\alpha_1 q}{y},y,z+\frac{\alpha_1 w}{y},w,q,\eta;\alpha_0+\alpha_1,-\alpha_1,\alpha_2+\alpha_1 \right),\\
s_2:(*) \rightarrow &\left(-x,y-\frac{2\alpha_2 q}{x}-\frac{q}{x^2},-z,-w+\frac{2\alpha_0}{x}+\frac{1}{x^2},-q,-\eta;\alpha_0,\alpha_1+2\alpha_2,-\alpha_2 \right).
\end{align*}
\end{theorem}

\begin{proposition}
Let us define the following translation operators{\rm : \rm}
\begin{align}\label{eq:13}
\begin{split}
&T_1:=s_1 s_2 s_1 s_0, \quad T_2:=s_1 T_1 s_1.
\end{split}
\end{align}
These translation operators act on parameters $\alpha_i$ as follows$:$
\begin{align}
\begin{split}
T_1(\alpha_0,\alpha_1,\alpha_2)=&(\alpha_0,\alpha_1,\alpha_2)+(-2,2,0),\\
T_2(\alpha_0,\alpha_1,\alpha_2)=&(\alpha_0,\alpha_1,\alpha_2)+(0,-2,2).
\end{split}
\end{align}
\end{proposition}

\section{Particular solution}
In this section, we study a solution of the system \eqref{eq:11} which is written by the use of known functions. 
\begin{proposition}
The system \eqref{eq:11} has the following invariant divisor\rm{:\rm}
\begin{center}\label{inv}
\begin{tabular}{|c|c|} \hline
parameter's relation  & invariant divisor \\ \hline
$\alpha_1=0$ & $y$ \\ \hline
\end{tabular}
\end{center}
\end{proposition}
Under the condition $\alpha_1=0$, elimination of $y$ from the system \eqref{eq:11} gives
\begin{equation}\label{eq:14}
  \left\{
  \begin{aligned}
\frac{dx}{dt}=&-(xw-\alpha_2)x+\frac{1}{2},\\
\frac{dz}{dt}=&-(zq-\alpha_0)z-\frac{\eta}{2},\\
\frac{dw}{dt}=&(xw-zq-\alpha_2)w,\\
\frac{dq}{dt}=&(zq-xw-\alpha_0)q.
   \end{aligned}
  \right. 
\end{equation}
At first, we find a particular solution $(x,z,w,q)=(x,z,0,0)$, and the system in $(x,z)$ are given by
\begin{equation}\label{eq:15}
  \left\{
  \begin{aligned}
\frac{dx}{dt}=&\alpha_2 x+\frac{1}{2},\\
\frac{dz}{dt}=&\alpha_0 z-\frac{\eta}{2}.\\
   \end{aligned}
  \right. 
\end{equation}
This system can be solved by
\begin{equation}\label{eq:16}
  \left\{
  \begin{aligned}
x(t)=&C_1 e^{\alpha_2 t}-\frac{1}{2\alpha_2},\\
z(t)=&C_2 e^{\alpha_0 t}+\frac{\eta}{2\alpha_0},
   \end{aligned}
  \right. 
\end{equation}
where $C_1,C_2$ are integral constants.

Next, we find a particular solution $(x,z,w,q)=(x,z,w,0)$, and the system in $(x,z,w)$ are given by
\begin{equation}\label{eq:17}
  \left\{
  \begin{aligned}
\frac{dx}{dt}=&-(xw-\alpha_2)x+\frac{1}{2},\\
\frac{dz}{dt}=&\alpha_0 z-\frac{\eta}{2},\\
\frac{dw}{dt}=&(xw-\alpha_2)w.
   \end{aligned}
  \right. 
\end{equation}
Elimination of $w$ from the system in the variables $x,w$ gives
\begin{equation}
\frac{d^2 x}{dt^2}=\frac{1}{x} \left(\frac{dx}{dt} \right)^2-\frac{\alpha_2}{2}-\frac{1}{4x}.
\end{equation}
This system can be solved by the following two solutions
\begin{equation}\label{eq:18}
x(t)=\frac{e^{-C_1(t+C_2)} \{ (e^{C_1(t+C_2)}-\alpha_2)^2-C_1^2 \}}{4C_1^2} 
\end{equation}
and
\begin{equation}\label{eq:19}
x(t)=\frac{e^{-C_1(t+C_2)} \{ e^{2C_1(t+C_2)}(\alpha_2^2-C_1^2)-2\alpha_2 e^{C_1(t+C_2)}+1 \}}{4C_1^2}, 
\end{equation}
where $C_1,C_2$ are integral constants. Of course, the system in the variable $z$ can be solved by \eqref{eq:16}.

\section{Polynomial Hamiltonian}

For the system \eqref{eq:11} let us try to seek its first integrals of polynomial type with respect to $x,y,z,w,q$.
\begin{proposition}
This system \eqref{eq:11} has its first integral\rm{:\rm}
\begin{equation*}
\frac{d(y-wq)}{dt}=-(y-wq).
\end{equation*}
\end{proposition}
By solving this equation, we can obtain
\begin{equation*}
y-wq=e^{-t}.
\end{equation*}
By using this, we show that elimination of $y$ from the system \eqref{eq:11} gives a polynomial Hamiltonian system.

\begin{theorem}
By the transformations
\begin{equation*}
q_1 =w, \quad p_1 =x, \quad q_2 =e^{t} q, \quad p_2 =\frac{z}{e^{t}}, \quad s=\frac{1}{e^{t}},
\end{equation*}
elimination of $y$ from the system \eqref{eq:11} gives a 2-parameter family of coupled Hamiltonian systems in dimension four explicitly given by
\begin{equation}\label{eq:3}
  \left\{
  \begin{aligned}
   \frac{dq_1}{ds} =&\frac{\partial H}{\partial p_1}=-\frac{q_1^2 p_1}{s}+\frac{\alpha_2 q_1}{s}-\frac{p_2}{s},\\
   \frac{dp_1}{ds} =&-\frac{\partial H}{\partial q_1}=\frac{q_1p_1^2}{s}-\frac{\alpha_2 p_1}{s}-\frac{1}{2s},\\
   \frac{dq_2}{ds} =&\frac{\partial H}{\partial p_2}=-\frac{q_2^2 p_2}{s}-\frac{(\alpha_1+\alpha_2)q_2}{s}-\frac{p_1}{s},\\
   \frac{dp_2}{ds} =&-\frac{\partial H}{\partial q_2}=\frac{q_2 p_2^2}{s}+\frac{(\alpha_1+\alpha_2)p_2}{s}+\frac{\eta}{2}
   \end{aligned}
  \right. 
\end{equation}
with the polynomial Hamiltonian
\begin{align}\label{eq:4}
\begin{split}
H=&K_1(q_1,p_1,s;\alpha_2)+K_2(q_2,p_2,s;\alpha_1+\alpha_2)-\frac{p_1p_2}{s}\\
=&-\frac{q_1^2 p_1^2-2\alpha_2 q_1p_1-q_1}{2s}-\frac{q_2^2 p_2^2+2(\alpha_1+\alpha_2)q_2p_2+\eta s q_2}{2s}-\frac{p_1p_2}{s}.
\end{split}
\end{align}
\end{theorem}
The symbols $K_1(q_1,p_1,s;\alpha)$ and $K_2(q_2,p_2,s;\alpha)$ denote
\begin{align}
\begin{split}
K_1(q_1,p_1,s;\alpha)=&-\frac{q_1^2 p_1^2-2\alpha q_1p_1-q_1}{2s},\\
K_2(q_2,p_2,s;\alpha)=&-\frac{q_2^2 p_2^2+2\alpha q_2p_2+\eta s q_2}{2s}.
\end{split}
\end{align}

This Hamiltonian system can be considered as a 1-parameter family of coupled polynomial Hamiltonian systems in dimension four.

We remark that for this system we tried to seek its first integrals of polynomial type with respect to $q_1,p_1,q_2,p_2$. However, we can not find. Of course, the Hamiltonian $H$ is not the first integral.

\begin{proposition}
The system \eqref{eq:3} is equivalent to the coupled equations:
\begin{equation}
  \left\{
  \begin{aligned}
   \frac{d^2y}{ds^2} =&\frac{1}{y}\left(\frac{dy}{ds} \right)^2-\frac{1}{s} \frac{dy}{ds}-\frac{\alpha_2}{2s^2}-\frac{1}{4s^2 y}-\frac{y^2 w}{s^2},\\
   \frac{d^2 w}{ds^2} =&\frac{1}{w}\left(\frac{dw}{ds} \right)^2-\frac{1}{s} \frac{dw}{ds}+\frac{\alpha_0 \eta}{2s}-\frac{\eta^2}{4w}-\frac{yw^2}{s^2},
   \end{aligned}
  \right. 
\end{equation}
where $y:=p_1$ and $w:=p_2$.
\end{proposition}

We study two Hamiltonians $K_1$ and $K_2$ in the principal parts of the Hamiltonian $H$.

At first, we study the Hamiltonian system
\begin{align}\label{44}
\begin{split}
\frac{dq_1}{ds}&=\frac{\partial K_1}{\partial p_1}, \quad \frac{dp_1}{ds}=-\frac{\partial K_1}{\partial q_1}
\end{split}
\end{align}
with the polynomial Hamiltonian
\begin{align}
\begin{split}
K_1:=-\frac{q_1^2 p_1^2-2\alpha q_1p_1-q_1}{2s},
\end{split}
\end{align}
where setting $q_2=p_2=0$ in the Hamiltonian $H$, we obtain $K_1$.

The system has the first integral $I_1${\rm : \rm}
\begin{equation}
I_1=q_1^2 p_1^2-2\alpha q_1p_1-q_1.
\end{equation}

Next, we study the Hamiltonian system
\begin{align}\label{444}
\begin{split}
\frac{dq_2}{ds}&=\frac{\partial K_2}{\partial p_2}, \quad \frac{dp_2}{ds}=-\frac{\partial K_2}{\partial q_2}
\end{split}
\end{align}
with the polynomial Hamiltonian
\begin{align}
\begin{split}
K_2(q_2,p_2,s;\alpha)=&-\frac{q_2^2 p_2^2+2\alpha q_2p_2+\eta s q_2}{2s},
\end{split}
\end{align}
where setting $q_1=p_1=0$ in the Hamiltonian $H$, we obtain $K_2$.

{\bf Step 1:} We make the change of variables:
\begin{equation}
x_1=s q_2, \quad y_1=\frac{p_2}{s}.
\end{equation}

Then, we can obtain the system with the polynomial Hamiltonian:
\begin{align}
\begin{split}
\tilde{K}_2(x_1,y_1,s;\alpha)=&-\frac{x_1^2 y_1^2+2(\alpha-1) x_1y_1+\eta x_1}{2s}.
\end{split}
\end{align}
This system has the first integral $I_2${\rm : \rm}
\begin{equation}
I_2=x_1^2 y_1^2+2(\alpha-1) x_1y_1+\eta x_1.
\end{equation}

We also study its symmetry.

\begin{theorem}\label{th:2}
The system \eqref{eq:3} admits extended affine Weyl group symmetry of type $D_3^{(2)}$ as the group of its B{\"a}cklund transformations, whose generators $s_0,s_1,s_2,\pi$ defined as follows$:$ with {\it the notation} $(*):=(q_1,p_1,q_2,p_2,\eta,s;\alpha_0,\alpha_1,\alpha_2):$
\begin{align}
\begin{split}
s_0:(*) \rightarrow &\left(q_1,p_1,q_2-\frac{2\alpha_0}{p_2}+\frac{\eta s}{p_2^2},p_2,\eta,-s;-\alpha_0,\alpha_1+2\alpha_0,\alpha_2 \right),\\
s_1:(*) \rightarrow &\left(q_1,p_1+\frac{\alpha_1 q_2}{q_1q_2+1},q_2,p_2+\frac{\alpha_1 q_1}{q_1q_2+1},\eta,s;\alpha_0+\alpha_1,-\alpha_1,\alpha_2+\alpha_1 \right),\\
s_2:(*) \rightarrow &\left(-q_1+\frac{2\alpha_2}{p_1}+\frac{1}{p_1^2},-p_1,-q_2,-p_2,-\eta,-s;\alpha_0,\alpha_1+2\alpha_2,-\alpha_2 \right),\\
\pi:(*) \rightarrow &\left(-\eta sq_2,-\frac{p_2}{\eta s},-\frac{q_1}{\eta s},-\eta s p_1,\eta,s;\alpha_2,\alpha_1,\alpha_0 \right),
\end{split}
\end{align}
where $\pi$ is its diagram automorphism of Dynkin diagram of type $D_3^{(2)}$.
\end{theorem}


\begin{thebibliography}{99}

\bibitem[1]{Sasano3} Y. Sasano, {\em Symmetries in the systems of types $D_3^{(2)}$ and $D_5^{(2)}$}, preprint.

\bibitem[2]{Cosgrove} C. M. Cosgrove,
{\em Higher order Painlev\'e equations in the polynomial class II, Bureau symbol P1}, Studies in Applied Mathematics. {\bf 116} (2006).

\bibitem[3]{8} P. Painlev\'e, {\em M\'emoire sur les \'equations diff\'erentielles dont l'int\'egrale g\'en\'erale est uniforme}, Bull. Soci\'et\'e Math\'ematique de France. {\bf 28} (1900),  201--261.

\bibitem[4]{9} P. Painlev\'e, {\em Sur les \'equations diff\'erentielles du second ordre et d'ordre sup\'erieur dont l'int\'egrale est uniforme}, Acta Math. {\bf 25} (1902), 1--85. 

\bibitem[5]{10} B. Gambier, {\em Sur les \'equations diff\'erentielles du second ordre et du premier degr\'e dont l'int\'egrale g\'en\'erale est \`a points critiques fixes}, Acta Math. {\bf 33} (1910), 1--55.


\bibitem[6]{Cosgrove2} C. M. Cosgrove,
{\em All binomial-type Painlev\'e equations of the second order and degree three or higher}, Studies in Applied Mathematics. {\bf 90} (1993), 119-187.

\bibitem[7]{13} F. Bureau, 
{\em Integration of some nonlinear systems of ordinary differential equations}, 
Annali di Matematica. {\bf 94} (1972), 345--359. 

\bibitem[8]{14} J. Chazy, 
{\em Sur les \'equations diff\'erentielles dont l'int\'egrale g\'en\'erale est uniforme et admet des singularit\'es essentielles mobiles}, 
Comptes Rendus de l'Acad\'emie des Sciences, Paris. {\bf 149} (1909), 563--565. 

\bibitem[9]{15} J. Chazy, 
{\em Sur les \'equations diff\'erentielles dont l'int\'egrale g\'en\'erale poss\'ede une coupure essentielle mobile }, 
Comptes Rendus de l'Acad\'emie des Sciences, Paris. {\bf 150} (1910), 456--458. 


\bibitem[10]{16} J. Chazy, 
{\em Sur les \'equations diff\'erentielles du trousi\'eme ordre et d'ordre sup\'erieur dont l'int\'egrale a ses points critiques fixes}, 
Acta Math. {\bf 34} (1911), 317--385.




\end{thebibliography}
\end{document}